\documentclass[11pt]{amsart}
 \usepackage{amsmath,amssymb,amsthm,amsfonts}
 \usepackage{epsfig} %% add this line and the next only if you have pictures
 \usepackage{graphics} %% pictures should be in esp format
\newcommand{\abs}[1]{| #1 |}
\newcommand{\Abs}[1]{\left| #1\right|}

\def\C {\mathbb C}
\newcommand{\cl} {\overline}
\newcommand{\dd}{\delta}

\newcommand{\D}{\mathbb D}
\newcommand{\e}{\epsilon}
\newcommand{\Frac}{\displaystyle\frac}

\newcommand{\Jac}{\operatorname{Jac}}

\newcommand{\norm} [1]{\left\| #1\right\|}
\newcommand{\of}{\circ}
\newcommand{\p}{\partial}
\newcommand{\R}{\mathbb R}
\newcommand{\re}{\text{\rm Re}\,}

\newcommand{\set}[1]{\left\{ #1\right\}}

\newcommand{\To}{\longrightarrow}

\newcommand{\W}{\Omega}
\newcommand{\z}{\zeta}

\newtheorem{theorem}{Theorem}

\newtheorem{lemma}{Lemma}
\newtheorem{prop}{Proposition}

\theoremstyle{definition}
\newtheorem{definition}{Definition}

\theoremstyle{remark}
\newtheorem{remark}{Remark}

\title[Kobayashi metric]
      {Asymptotic Behavior of the Kobayashi Metric on Convex Domains}

\author{Lina Lee}

\email{ linalee@umich.edu}

\begin{document}

\maketitle

\begin{abstract}
In this paper, we calculate estimates for invariant metrics on a finite type convex domain in $\C^n$ using the Sibony metric. We also discuss a possible modification of the Sibony metric.
  \end{abstract}
\section{Introduction}

The Kobayashi metric $F(P,\xi)$ on a domain $\W\subset\C^n$ at a point $P\in\W$ in the direction $\xi\in T_P(\W)$ is defined as follows:
\begin{equation}\label{1157}
F(P,\xi)=\inf\set{\alpha>0:\exists \phi\in\W(\D),\;\phi(0)=P,\;\phi'(0)=\xi/\alpha},
\end{equation}
where $\W(\D)$ denotes the family of holomorphic mappings from the unit disc $\D$ in $\C$ to $\W$. It is known that the Kobayashi metric is greater than any biholomorphically invariant metric $G$ that satisfies the following properties:
\begin{enumerate}
\item
$G^\D:\D\times\C\To \R^+\cup\set 0$ coincides with the Poincar\'e metric on the unit disc in $\C$.
\item
$G$ is non-increasing under holomorphic mappings, i.e., if $\Phi:\W\To\tilde\W$ is a holomorphic mapping and $P\in\W$, $\xi\in T_P(\W)$, then 
$$
G^\W(P,\xi)\ge G^{\tilde\W}(\Phi(P),\Phi_*(P)\xi).
$$
\end{enumerate}

It has been of importance to study the asymptotic behavior of the Kobayashi metric near the boundary of a holomorphically convex domain. Several authors have proved results on pseudoconvex domains: Ian Graham proved a result on a strongly pseudoconvex domain \cite{Graham} and David Catlin studied the behavior on a weakly pseudocovnex domain in $\C^2$ \cite{Catlin}. 

As we can see from the definition of the Kobayashi metric (\ref{1157}), the difficulty in estimating the Kobayashi metric lies in finding the lower estimate since the upper estimate can be found rather easily by constructing one analytic disc in $\W$ that satisfies the desired properties. 

Graham \cite{Graham} calculated the metric explicitly on ellipsoids and found the estimate on a strongly pseudoconvex domain by approximating it with ellipsoids and proving results on localization of the metric. Catlin \cite{Catlin} proved the result by estimating Carath\'eodory metric $F_C(P,\xi)$, which is defined as follows:
$$
F_C(P,\xi)=\sup\set{ \abs{f_*(P)\xi}=\Abs{\sum_{i=1}^n\frac{\p f(P)}{\p z_i}\xi_i}: f\in\D(\W),\; f(P)=0}.
$$

The Carath\'eodory metric satisfies above two properties and hence is less than the Kobayashi metric. So one can estimate the Carath\'eodory metric and find a lower estimate for the Kobayashi metric. 

In this paper, we estimate the Kobayashi metric on a convex domain in $\C^n$ using the Sibony metric, whose definition can be found in section 2. The Sibony metric also satisfies the two properties above and hence gives a lower estimate for the Kobayashi metric. The advantage of using the Sibony metric over using the Carath\'eodory metric is that the Sibony metric uses bounded plurisubharmonic functions whereas the Carath\'eodory metric uses bounded holomorphic functions, which are usually more difficult to construct than plurisubharmonic functions. In section 2, we give a more detailed explanation of the Sibony metric. 
 
We assume $\W=\set{\rho<0}\subset\subset\C^n$ is a smoothly bounded convex domain, $P\in\p\W$, $\nu=\nabla\rho(P)/\norm{\nabla\rho(P)}$ and let $P_\dd=P-\dd\nu\in\W$. 

For $\xi\in T_P^\C(\p\W)= T_p(\p\W)\cap JT_P(\p\W)$, where $J$ is the standard complex structure of $\C^n$, we define $\Delta (\p\W,P,\xi)$ as the tangency of the $\p\W$ at $P$ in the direction $\xi$, i.e.,
\begin{equation}\label{337}
\Delta(\p\W,P,\xi)=v_0(\rho(P+\xi z)), \quad z\in\C
\end{equation}
where $v_0(f(z))$ denotes the vanishing order of $f$ at $z=0$.

L. Lempert proved that the Kobayashi metric and the Carath\'eodory metric coincide on a convex domain in $\C^n$ \cite{Lempert}. Hence the Sibony metric also coincides with the Carath\'eodory metric and the Kobayashi metric. Let us denote the (Kobayashi or Sibony or Carath\'eodory) metric as $F(Q,\xi)$ for $Q\in\W$ and $\xi\in T_Q(\W)$. Then we have the following theorem. 
\begin{theorem}\label{614}
If $\W\subset\C^n$ is a smoothly bounded convex domain of finite type, then we have
\begin{equation}\label{1200}
F(P_\dd,\xi)\approx\frac{\abs\xi}{\dd^{1/m}},\quad\xi\in T_P^\C(\p\W),
\end{equation}
where $m=\Delta(\p\W,P,\xi)$, and
\begin{equation}\label{1201}
F(P_\dd,\nu)\approx\frac{1}{\dd}
\end{equation}
for all sufficiently small $\dd>0$.
\end{theorem}
\begin{remark}
In (\ref{1200}) and (\ref{1201}), the notation ``$\approx$'' means that there exist positive constants $c$, $C$, $c'$ and $C'$ that do not depend on $\dd$ such that
\begin{gather*}
c\frac{\abs\xi}{\dd^{1/m}}\le F(P_\dd,\xi)\le C\frac{\abs\xi}{\dd^{1/m}},\quad{and}\\
c'\frac{1}{\dd}\le F(P_\dd,\nu)\le C' \frac{1}{\dd},
\end{gather*}
for all $\dd>0$ sufficiently small.
\end{remark}
The boundedness from above can be easily shown using the definition of the Kobayashi metric. We can express the defining function using the Taylor series and find an analytic disc that has the proper size in the estimating direction. For more details, refer \cite{Lina}. In section 2, we prove the boundedness from below. 

We also prove the following theorem. 
 
\begin{theorem}\label{615}
Suppose $\W\subset\subset\C^n$ is a smoothly bounded convex domain. Let $X=a\nu+bT$, where $T\in T_P(\W)$ and $a,b>0$. Then we have  
$$
F(P_\dd, X)\ge \frac{\abs a}{6\dd}.
$$
\end{theorem}

In section 2, we give a brief background of invariant metrics and finite type and, in section 3, we prove Theorem \ref{614} and Theorem \ref{615}. In section 4, we discuss a possible modification of the Sibony metric. 

\section{Background: Invariant Metrics and the Concept of Finite Type}

We say $F: T\W\To\R^+\cup\set{0}$ is an invariant metric if $F$ is invariant under biholomorphic mappings, i.e., if $\Phi:\W_1\To\W_2$ is a biholomorphic mapping between $\W_1$ and $\W_2$ and $P\in\W_1$, $\xi\in T_P(\W_1)$, then
\begin{equation}\label{1256}
F^{\W_1}(P,\xi)=F^{\W_2}(\Phi(P),\Phi_*(P)\xi).
\end{equation}

For example, the Poincar\'e metric $P(z,\xi)$ on the unit disc $\D$ in $\C$, which is defined as
$$
P(z,\xi)=\frac{\abs\xi}{1-\abs z^2},
$$
is invariant under automorphisms of the unit disc.

Two possible generalizations of the Poincar\'e metric to an arbitrary domain $\W$ in $\C^n$ are the Kobayashi metric, $F_K(P,\xi)$, and the Carath\'edory metric, $F_C(P,\xi)$, which are defined as follows:
\begin{gather}
F_K(P,\xi)=\inf\set{\alpha: \exists \phi\in\W(\D), \;\phi(0)=P,\; \phi'(0)=\xi/\alpha,\;\alpha>0};\\
F_C(P,\xi)=\sup\set{ \abs{f_*(P)\xi}=\Abs{\sum_{i=1}^n\frac{\p f(P)}{\p z_i}\xi_i}: f\in\D(\W),\; f(P)=0},
\end{gather}
where $A(B)$ denotes the family of holomorphic mappings from $B$ to $A$ and $\D$ the unit disc in $\C$.

The Kobayashi metric is the largest pseudometric and the Carath\'eodory metric is the smallest in the following sense:
\begin{prop}
Suppose that $\tilde F^\W:T\W\To\R^+\cup\set 0$ is a pseudometric on $\W$ such that $\tilde F^\D$ coincides with the Poincar\'e metric  and $\tilde F$ is non-increasing under holomorphic mappings, i.e., if $\Phi:\W_1\To\W_2$ is a holomorphic mapping and $P\in \W_1$, then we have
$$
\tilde F^{\W_1}(P,\xi)\ge \tilde F^{\W_2}(\Phi(P),\Phi_*(P)\xi),\quad \forall\xi\in T_P^\C(\W_1).
$$
Then we always have $F_C^\W(P,\xi)\le \tilde F^\W(P,\xi)\le F_K^\W(P,\xi)$.
\end{prop}
The Sibony metric is defined as follows:
\begin{definition}[Sibony metric]\label{Sibony-definition}
Let $\W\in\C^n$ be a domain and $P\in\W$. We define a set of
functions, $A_\W(P)$,  such that $u\in A_\W(P)$ if and only if
\begin{enumerate}
\item
$u$ is $C^2$ near $P$;
\item
$u(P)=0$;
\item
$ 0\le u(z)\le 1$ for all $z\in\W$;
\item
$\log u$ is plurisubharmonic on $\W$.
\end{enumerate}
We define the infinitesimal Sibony metric $F_\W^S$ at $P$ in the
direction $\xi\in\C^n$ as follows:
\begin{equation}
F_S(P,\xi)\equiv \sup_{u\in
A_\W(P)}\left(\sum_{i,j=1}^n\frac{\p^2u}{\p z_i\p \cl
z_j}(P)\xi_i\cl\xi_j\right)^{\frac{1}{2}}.\label{Sibony-metric}
\end{equation}
\end{definition}
The Sibony metric coincides with the Poincar\'e metric on the unit disc and is non-increasing under holomorphic mappings. Hence we have that
$$
F_C^\W(P,\xi)\le F_S^\W(P,\xi)\le F_K^\W(P,\xi).
$$
\paragraph{Finite Type}
We say the boundary of a domain in $\C^n$ if of finite type if the maximum tangency of the boundary with any one dimensional holomorphic variety is finite, i.e., 
$$
\sup\set{\frac{v_0(\rho\of\phi)}{v_0(\phi)}:\phi\in\C^n(\D),\;\phi(0)=P}<\infty.
$$
For more details, refer \cite{KrantzSCV} and \cite{D'Angelo}. 
McNeal \cite{McNeal} showed that the finite type condition of a boundary of a convex domain in $\C^n$ is same as the finite type condition with $\phi$ replaced with complex lines through $P$. So we introduced the notation $\Delta(\p\W,P,\xi)$ in (\ref{337}), which actually gives you the type in the direction $\xi$ if $\W$ is convex. 

\section{Estimation on a Convex Domain}
Throughout this section we assume that $\W=\set{\rho<0}\subset\subset\C^n$ is a smoothly bounded convex domain, $P\in\p\W$, $\xi\in T_P^\C(\p\W)$ and $\nu=\nabla\rho(P)/\norm{\nabla\rho(P)}$, which is the outwarad unit normal vector at $P$. Let $P_\dd=P-\dd\nu$.

We use the following lemma by Bruna, Nagel and Wainger proved in \cite{BNW}. 

\begin{lemma}[Bruna, Nagel, Wainger]\label{BNWlemma}
Let us define a set of functions on $\R$ as follows:
\begin{equation}
C(m,r)\equiv\set{f(x)=a_2x^2+\cdots a_m x^m: a_i\in\R, f''(x)\ge 0 \;\forall x\in [0,r]}. \label{BNW}
\end{equation}
Then there exists a constant $C$ such that
$$
f(x)\ge C(\abs{a_2} x^2+\cdots+\abs{a_m}x^m),\quad\forall f\in C(m,r),\;\forall x\in[0,r].
$$
\end{lemma}
\begin{prop}Let $\W=\set{\rho<0}\subset\subset\C^n$ be a smoothly bounded convex domain, $P\in\p\W$ and $\nu$ the unit outward real normal vector to  $\p\W$ at $P$ with $\norm\nu=1$. Let $\xi\in T_P^\C(\p\W)$, $\norm\xi=1$, and $\Delta(P,\p\W,\xi)>2$. If we let
$$
R_\xi(\dd):=\sup\set{\abs z : P-\dd\nu +\xi z\in\W,\, z\in\C},
$$ then
$$
R_\xi(\dd)\approx \dd^{1/m}
$$
for $\dd>0$ sufficiently small.
\end{prop}
\begin{proof}
We may assume $P=0$, $\nabla\rho(P)=(0,\dots, 1)$ and $\xi=(1,0,\dots,0)$. Then near $P=0$, $\rho$ can be expressed as  $\rho=\re z_n+O(\abs z^2)$ and, if we evaluate $\rho$ at $(0,\dots,0,-\dd)$ in the $z_1$-direction, we get
$$
\rho((\z,0,\dots,0,-\dd))=-\dd+\sum_{p+q_1+q_2=2,\, p\ge 1}^{m-1}a_{pq_1q_2}\dd^p\z^{q_1}\cl\z^{q_2}+O((\dd^2+\abs\z^2)^{\frac{m}{2}}).
$$
Let $\abs\z=c\dd^{1/m}$.
Then
\begin{multline}
\rho((\z,0,\dots,0,-\dd))=-\dd+\sum_{p+q=2,\,p\ge 1}^{m-1}b_{pq}\dd^p(c\dd^{\frac{1}{m}})^q+O((\dd^2+(c\dd^{1/m})^2)^{\frac{m}{2}})\\
\le -\dd+\sum_{p+q=2,\,p\ge 1}^{m-1}b_{pq}\dd^p(c\dd^{\frac{1}{m}})^q+C(\dd^m+c'\dd)<0,
\end{multline}
for some constants $C$ and $c'$. Hence, we get
$$
R_\xi(\dd)\gtrsim\dd^{1/m}.
$$
Next, we will show that, for a fixed $\e\in (0,1/m)$,  there does not exist a constant $c$ such that, for all  sufficiently small $\dd$,
$$
\rho((\z,0,\dots,0,-\dd))<0,\quad\abs\z=c\dd^{\frac{1}{m}-\e}.
$$
Let $\abs\z=c\dd^{1/m-\e}$ and look at the Taylor expansion:
\begin{multline}
\rho((\z,0,\dots,0,-\dd))\\
=-\dd+\sum_{p+q=2,\,p\ge 1}^{m} b_{pq}\dd^p(c\dd^{\frac{1}{m}-\e})^q+b_{0m}(c\dd^{\frac{1}{m}-\e})^m+O((\dd^2+(c\dd^{\frac{1}{m}-\e})^2)^{\frac{m+1}{2}}).
\end{multline}
By Lemma \ref{BNWlemma}, we get
\begin{multline}
\rho((\z,0,\dots,0,-\dd))\\
\ge -\dd+C\left(\sum_{p+q=2,\,p\ge 1}^m\abs{b_{pq}}\dd^p(c\dd^{\frac{1}{m}-\e})^q+\abs{b_{0m}}(c\dd^{\frac{1}{m}-\e})^m\right)\\
-C'(\dd^{m+1}+c'(\dd^{\frac{1}{m}-\e})^{m+1})\\
\ge -\dd+C\abs{b_{0m}}(c\dd^{\frac{1}{m}-\e})^m-C'(\dd^{m+1}+c'(\dd^{\frac{1}{m}-\e})^{m+1})\\
=C''\dd^{1-\e m}-\dd-C'(\dd^{m+1}+c'(\dd^{\frac{1}{m}-\e})^{m+1}>0,
\end{multline}
for $\dd$ sufficiently small.
Hence, we conclude that
$$
R_\xi(\dd)\approx\dd^{1/m}.
$$
\end{proof}
To find the lower bound for the Sibony metric, we construct a plurisubharmonic function that satisfies the conditions of the definition and has a large Hessian in the $\xi$-direction. We could find such a plurisubharmonic function by modifying the construction of a plurisubharmonic function in \cite{McNeal}.
\begin{prop}\label{943}
Suppose $\xi\in T_P^\C(\p\W)$ and $\Delta(\p\W,P,\xi)=m$. Then
$$
F(P_\dd,\xi)\gtrsim \frac{\abs\xi}{\dd^{1/m}}.
$$
\end{prop}
\begin{proof}
We may assume $\xi$ is the $z_1$ direction. We let the $\re z_1$-direction be such that the distance from $P-\dd\nu$ to the boundary along $\re z_1$ axis will be the greatest among all distances between $P-\dd\nu$ and the boundary along the $z_1$ axis, i.e.,
\begin{multline}
\sup\set{r>0:\rho((r,0,\dots,0,-\dd))\in\W}\\
=\sup\set{r>0:\rho((re^{i\theta},0,\dots,0,-\dd))\in\W,\;0\le \theta<2\pi}.
\end{multline}
Let $R$ be such a  distance
$$
R=\sup\set{r>0:\rho((r,0,\dots,0,-\dd))\in\W}.
$$
Then by the Proposition, we know that
$$
R\approx \dd^{1/m}.
$$
Let $Q=(R,0,\dots,0,-\dd)$. Now we will show that
$$
\frac{\p\rho}{\p z_1}(Q)\approx\dd^{1-1/m}.
$$
using the technique in \cite{McNeal}.

Consider the real tangent space to $\p\W$ at $Q$:
\begin{multline}
\re\Big(\frac{\p\rho}{\p z_1}(Q)(z_1-R)+\frac{\p\rho}{\p z_2}(Q) z_2+\cdots\\
+\frac{\p\rho}{\p z_{n-1}}(Q)z_{n-1}+\frac{\p\rho}{\p z_n}(Q)(z_n+\dd)\Big)=0.
\end{multline}
Let $S$ be the intersection point between the above tangent space and the $\re z_n$ axis, i.e.,
\begin{multline}
S=\Bigg\{z:\re\Bigg(\frac{\p\rho}{\p z_1}(Q)(z_1-R)+\frac{\p\rho}{\p z_2}(Q) z_2+\cdots\\
+\frac{\p\rho}{\p z_{n-1}}(Q)z_{n-1}+\frac{\p\rho}{\p z_n}(Q)(z_n+\dd)\Bigg)=0\Bigg\}
\cap\set{(0,\dots,0,x),\;x\in\R}.
\end{multline}
If we let $S=(s,0,\dots,0)$, then, by convexity, we have
$$
\abs{S-P}=\abs{(s,0,\dots,0)-0}= s\ge 0.
$$
Therefore, if we evaluate the tangent space at $S$, we get
$$
\re\left(\frac{\p\rho}{\p z_1}(Q)(-R)+\frac{\p\rho}{\p z_n}(Q)(s+\dd)\right)=0
$$
Hence we have
$$
\Abs{\frac{\p\rho}{\p z_1}(Q)}R\ge\Abs{\re\frac{\p\rho}{\p z_1}(Q)R}=\Abs{\frac{\p\rho}{\p z_n}(Q)}(s+\dd)\ge \Abs{\frac{\p\rho}{\p z_n}(Q)}\dd\approx\dd.
$$
\begin{equation}\label{pshfunction1}
\Abs{\frac{\p\rho}{\p z_1}(Q)}\gtrsim\frac{\dd}{R}\approx\frac{\dd}{\dd^{1/m}}=\dd^{1-1/m}.
\end{equation}

To show the other direction, we look at the Taylor expansion of $\rho$ at $0$ and evaluate it along $z_1$-direction.
\begin{multline}
\rho((z,0,\dots,-\dd))=-\dd+\sum_{p+q_1+q_2=2,\; p\ge 1}^m a_{p q_1 q_2}\dd^p z^{q_1}\cl z^{q_2}\\
+\sum_{r_1+r_2=m} a_{r_1 r_2}z^{r_1}\cl z^{r_2}+O((\dd^2+\abs z^2)^{\frac{m+1}{2}}).
\end{multline}
Differentiating along the $z_1$ direction, we get
$$
\frac{\p\rho}{\p z_1}(Q)=\sum_{p+q=2,p\ge 1}^m b_{pq}\dd^p R^{q-1}+ b_m R^{m-1}+O((\dd^2+\abs z^2)^{\frac{m}{2}}).
$$
Since $R\le C\dd^{1/m}$, we get
\begin{multline}\label{pshfunction2}
\Abs{\frac{\p\rho}{\p z_1}(Q)}\lesssim\sum_{p+q=2,\,p\ge 1}^m\abs{b_{pq}}C^{q-1}\dd^{p+\frac{q-1}{m}}+\abs{b_m}C^{m-1}\dd^{1-1/m}+O(\dd^m+C^m\dd)\\
\lesssim\dd^{1-1/m}.
\end{multline}
By (\ref{pshfunction1}) and (\ref{pshfunction2}), we get
$$
\Abs{\frac{\p\rho}{\p z_1}(Q)}\approx\dd^{1-1/m}.
$$
Now we construct a candidate plurisubharmonic function for the Sibony metric $F_S(P-\dd\nu,\xi)$. Let
$$
f=\frac{1}{\dd}\left(\frac{\p\rho}{\p z_1}(Q) z_1+\cdots+\frac{\p\rho}{\p z_{n-1}}(Q) z_{n-1}+\frac{\p\rho}{\p z_n}(Q)(z_n+\dd)\right)
$$
and define $F_N$:
$$
F_N=f+\frac{1}{2!}f^2+\cdots+\frac{1}{N!}f^N,
$$
where the number $N$ will be chosen later.
If we let
$$
G_N=\abs{F_N}^2,
$$
then $G_N$ satisfies the following properties: $G_N(P-\dd\nu)=0$, $\log G_N$ is plurisubharmonic on $\W$, and
$$
\frac{\p^2 G_N}{\p z_1\p\cl z_1}(P-\dd\nu)=\frac{1}{\dd^2}\Abs{\frac{\p\rho}{\p z_1}(Q)}^2\approx\frac{\dd^{2-2/m}}{\dd^2}=\left(\frac{1}{\dd^{1/m}}\right)^2.
$$
Hence, if we can show that $G_N$ is bounded on $\W$, then we can conclude that
$$
F_S^\W(P_\dd,\xi)\gtrsim \frac{\abs\xi}{\dd^{1/m}}.
$$
Now we prove that $G_N$ is bounded and has an upper bound independent of $\dd$.

Since
$$
1+f+\frac{1}{2!}f^2+\cdots+\frac{1}{k!}f^k+\cdots=\exp f,
$$
we may find $N$ such that
$$
\abs{1+F_N(z)-\exp f(z)}<1,\quad\forall z\in\W.
$$
Therefore
$$
\abs{F_N}<1+\abs{\exp f -1}\le 2+e^{\re f}.
$$
Since $\re f=0$ defines a hyperplane and $\re f$ changes sign at $\re f=0$, we may assume that $\re f>0$ near the boundary and negative elsewhere. Then $e^{\re f}\le 1$ outside a small neighborhood of $P$. Therefore we have
$$
\re f\le \re f(Q)=\frac{1}{\dd}\re \left(\frac{\p\rho}{\p z_1}(Q)R\right)\lesssim\frac{1}{\dd}\dd^{1-1/m}\dd^{1/m}=1.
$$
Hence $G_N$ is uniformly bounded for all $\delta$.
\end{proof}

\begin{prop}\label{944}
$$
F(P_\dd,\nu)\ge\frac{1}{6\dd},
$$
\end{prop}
\begin{proof}
Let
$$
\rho(z)=2\re z_n+O(\abs z^2).
$$
Since $\W$ is convex, we see that  $\W\subset\set{\re z_n<0}$. Let us look at the function:
\begin{equation}\label{416}
u(z)=\frac{1}{9}\Abs{\frac{z_n+\dd}{z_n-\dd}}^2.
\end{equation}
Since $\re z_n<0$ for all $z\in\W$, the function inside the absolute value sign is holomorphic on $\W$. Hence $\log u$ is plurisubharmonic on $\W$. And $u(P-\dd\nu)=0$. We also have
$$
\Abs{\frac{z_n+\dd}{z_n-\dd}}\le 1+\Abs{\frac{2\dd}{z_n-\dd}}\le 1+\frac{2\dd}{\dd}=3,\quad z\in\W,
$$
since
$$
\abs{z_n-\dd}\ge\abs{\re z_n-\dd}=\abs{\re z_n}+\dd\ge\dd,\quad \forall z\in\W.
$$
Hence $0\le u(z)\le 1$ on $\W$. Finally,
$$
\left(\frac{\p^2 u(P_\dd)}{\p z_n\p\cl z_n}\right)^{1/2}=\left(\frac{1}{9}\frac{1}{4\dd^2}\right)^{1/2}=\frac{1}{6\dd}
$$
\end{proof}
Hence Theorem \ref{614} is proved by Proposition \ref{943} and Proposition \ref{944}.
%\begin{theorem}
%Let $X=a\nu+bT$, whre $T\in T_P(\p\W)$ and $a,b>0$. Then we have
%$$
%F(P_\dd, X)\ge \frac{\abs a}{6\dd}
%$$
%\end{theorem}

\bigskip
\noindent{\bf Proof of Theorem \ref{615}}
\begin{proof}
We use the same function $u(z)$ as in (\ref{416}):
$$
u(z)=\frac{1}{9}\Abs{\frac{z_n+\dd}{z_n-\dd}}^2.
$$
Then we have that
$$
F_S(P_\dd,X)\ge(\p\cl\p u(X,\cl X))^{1/2}=\frac{\abs a}{6\dd}.
$$
\end{proof}
\section{A Modification of the Sibony metric}

In this section, we discuss a possible modification of the Sibony metric.

\begin{definition}[Plurisubharmonic Metric]\label{psh-metric-def}
Let $\W\subset\C^n$ be a domain and $P\in\W$, $\xi\in\C^n$. We
define a set of functions $B_\W(P,\xi)$ such that $u\in B_\W(P,\xi)$
if and only if
\begin{enumerate}
\item
$u$ is $C^2$ near $P$;
\item
$u(P)=0$;
\item there exists a holomorphic disc $f:\D\To\cl\W$ such that $f(0)=P$, $\displaystyle f'(0)=\frac{\xi}{F_K^\W(P,\xi)}$ and $u$ satisfies
\begin{enumerate}
\item
$0\le u\of f(z)\le 1$, for all $z\in\D$;
\item
$\displaystyle \frac{u\of f(z)}{\abs z^2}$ is subharmonic on $\D$.
\end{enumerate}
\end{enumerate}
We define the plurisubharmonic metric $F_\W^P$ at $P\in\W$ in the
direction $\xi\in\C^n$ as follows:
\begin{equation}
F_P^\W(P,\xi)\equiv\sup_{u\in
B_\W(P,\xi)}\left(\sum_{i,j=1}^n\frac{\p^2 u}{\p z_i\p\cl
z_j}(P)\xi_i\cl\xi_j\right)^{\frac{1}{2}}.\label{psh-metric}
\end{equation}
\end{definition}

\begin{prop} If $\W\subset\subset\C^n$ is a pseudoconvex domain and $P\in\W$ and $\xi\in \C^n$, then
$$
F_S^\W(P,\xi)\le F_P^\W(P,\xi).
$$
\end{prop}
\begin{proof}
It is enough to show that the collection of candidate functions for the Sibony metric is a subset of collection of candidate functions for the plurisubharmonic metric, i.e., $A_\W(P)\subset B_\W(P,\xi)$. 

If $u\in A_\W(P)$, then, as stated in Definition
\ref{Sibony-definition}, $u$ is $C^2$ near $P$, $u(P)=0$, $0\le
u(z)\le 1$ for all $z\in\W$ and $\log u$ is plurisubharmonic on
$\W$. We need to show that $u$ satisfies the conditions of
Definition \ref{psh-metric-def}.

Since $\W$ is pseudoconvex, we can find an extremal disc
$f\in\W(\D)$ such that $f(0)=P$ and $\displaystyle
f'(0)=\frac{\xi}{F_K^\W(P,\xi)}$. Hence, since $0\le u\le 1$ on
$\W$, we get $0\le u\of f(z)\le 1$ for all $z\in\D$. We know that
$\log u$ is plurisubharmonic on $\W$. Therefore, $\log u\of f$ is
subharmonic on $\D$. Since $\log \abs z$ is harmonic, $\log u\of
f-\log\abs z^2$ is subharmonic on $\D$. Taking the exponential, we
see that  $\displaystyle \frac{u\of f(z)}{\abs z^2}$ is subharmonic
on $\D$.
\end{proof}

Next we will show that the metric $F_P^\W$ is invariant under
biholomorphic mappings. This is connected to the fact that the
Kobayashi metric is invariant under biholomorphic mappings.
\begin{prop}\label{psh-metric-invariant}
The plurisubharmonic metric is invariant under biholomorphic
mappings.
\end{prop}
\begin{proof}
Let $\W_1,\W_2\subset\C^n$, $P\in\W_1$ and $\xi\in\C^n$. Suppose that
$\Phi:\W_1\To\W_2$ is a biholomorphic mapping. We will show that
$$
u\of\Phi^{-1}\in B_{\W_2}(\Phi(P),\Phi_*(P)\xi), \quad \forall u\in
B_{\W_1}(P,\xi),
$$
where $B_\W(P,\xi)$ is the set of functions that satisfy the
conditions of Definition \ref{psh-metric-def}. Since $\Phi$ is a
biholomorphic mapping, $u\of\Phi^{-1}$ is $C^2$ near $\Phi(P)$ and
$u\of\Phi^{-1}(\Phi(P))=u(P)=0$. Thus the first two conditions are
satisfied. Now suppose that $f:\D\To\cl\W_1$ is a holomorphic curve
that satisfies the conditions of  Definition \ref{psh-metric-def}
for $F_P^{\W_1}(P,\xi)$. If we let $g(z)=\Phi\of f(z)$, then
$g(0)=\Phi(f(0))=\Phi(P)$ and
\begin{multline}
g'(0)=\Jac\Phi(P)f'(0)\\
=\Jac\Phi(P)\frac{\xi}{F_{\W_1}^K(P,\xi)}=\frac{\Phi_*(P)\xi}{F_{\W_1}^K(P,\xi)}=\frac{\Phi_*(P)\xi}{F_{\W_2}^K(\Phi(P),\Phi_*(P)\xi)}.
\end{multline}
The last equality holds since the Kobayashi metric is invariant
under $\Phi$.

Now we want to show that $u\of \Phi^{-1}$ satisfies the conditions
on the holomorphic curve $\Phi\of f$. But this is rather
straightforward since $f$ was chosen to be the holomorphic curve on
which the function $u$ satisfies the conditions of Definition
\ref{psh-metric-def} and $u\of\Phi^{-1}\of \Phi\of f(z)=u \of f(z)$.
\end{proof}

Now we want to prove that the plurisubharmonic metric coincides with
the Poincar\'e metric on the unit disc in $\C$.

\begin{prop} The plurisubharmonic metric coincides with the Poincar\'e metric on the unit disc in $\C$, i.e.,
$$
F_P(0,\xi)=P_\D(0,\xi)=\abs\xi.
$$
\end{prop}
\begin{proof}
By the following lemma, we know that $F_P(0,\xi)\le\abs\xi$. Also,
since $u(z)=\abs z^2$ satisfies the conditions of being a candidate
function, we get $F_P(0,\xi)=\abs\xi$.
\end{proof}
\begin{lemma}
If $u:\D\To\R$ satisfies
\begin{enumerate}
\item $u(0)=0$, $u$ is $C^2$ near $0$;
\item $0\le u\le 1$ on $\D$;
\item $\Frac{u(z)}{\abs z^2}$ is subharmonic on $\D$.
\end{enumerate}
then
$$
\frac{\p^2 u(0)}{\p z\p\cl z}\le 1.
$$
\end{lemma}
\begin{proof}
Since $u$ is $C^2$ near $0$ and has minimum at $0$, the first order
derivatives of $u$ at $0$ is $0$. Hence the Taylor expansion near
$0$ becomes
$$
u(z)=a\abs z^2+\re bz^2+O(\abs z^3),\quad a\in\R,\;b\in\C.
$$
Let $z=\abs z e^{i\theta}$. Then
$$
\frac{u(z)}{\abs z^2}=a+\re(b e^{i\theta})+O(\abs z).
$$
Since $\Frac{u(z)}{\abs z^2}$ is subharmonic on $\D$, by the  maximum
principle we get
$$
\frac{u(z)}{\abs z^2}\Big|_{z=0}\le\frac{u(z)}{\abs z^2}\Big|_{\abs
z=1}=u(z)|_{\abs z=1}\le 1.
$$
Therefore
$$
\frac{u(z)}{\abs z^2}\Big|_{z=0}=a+\re(b e^{2i\theta})\le 1.
$$
Choose $\theta_0$ such that $\re(be^{2i\theta_0})\ge 0$. We get
$$
a\le 1-\re(be^{2i\theta_0})\le 1.
$$
\end{proof}

\begin{prop}
The plurisubharmonic metric is less than or equal to the Kobayashi
metric, i.e.,
$$
F_P^\W(P,\xi)\le F_K^\W(P,\xi).
$$
\end{prop}
\begin{proof}
Let $u$ be the candidate function and $f$ be the curve that has the
derivative at the base point of the same size as the inverse of the
Kobayashi metric in the $\xi$ direction such that $0\le u\of f \le
1$ on $\D$ and $\Frac{u\of f(z)}{\abs z^2}$ is subharmonic on $\D$.
Then by the previous lemma, we get
$$
\frac{\p^2 u\of f(0)}{\p z\p\cl z}\le 1
$$
Rewriting the left hand side of the above inequality, we get
$$
\frac{\p^2 u\of f(0)}{\p z\p\cl z}=\sum_{j,k}\frac{\p^2
u(P)}{\p\xi_j\p\cl\xi_k}\Big(f'(0)\Big)_j\Big(cl{f'(0)}\Big)_k
=\frac{1}{(F_K(P,\xi))^2}\sum_{j,k}\frac{\p^2
u(P)}{\p\xi_j\p\cl\xi_k}\xi_j\cl\xi_k.
$$
Hence
$$
\left(\sum{j,k}\frac{\p^2 u(P)}{\p\xi_j\p\cl
\xi_k}\xi_j\cl\xi_k\right)^{1/2}\le F_K(P,\xi).
$$
Therefore, $F_P(P,\xi)\le F_K(P,\xi)$.
\end{proof}

\smallskip
\noindent Lina Lee\\
Mathematics Department\\
The University of Michigan\\
East Hall, Ann Arbor, MI 48109\\
USA\\
linalee@umich.edu\\

\end{document}